\documentclass{gtart}

\def\ifplaintex{\expandafter\ifx\csname documentclass\endcsname\relax}

\def\gtp{{\mathsurround=0pt\it $\cal G\mskip-2mu$eometry \&\ 
$\cal T\!\!$opology $\cal P\!$ublications}}  

\def\Addressesr{\bigskip
{\small \parskip 0pt \leftskip 0pt \rightskip 0pt plus 1fil \def\\{\par}
\sl\theaddress\par
\medskip
\rm Email:\stdspace\tt\theemail\hfill\rm Received:\qua\receiveddate \par}}

\def\recd{{\small Received:\qua\receiveddate\ifx\reviseddate\relax
\else\qquad Revised:\qua\reviseddate\fi\par}} 

\def\lognumber#1{\def\thelognumber{#1}}
\def\volumenumber#1{\def\thevolumenumber{#1}}
\def\volumeyear#1{\def\thevolumeyear{#1}}
\def\papernumber#1{\def\thepapernumber{#1}}
\def\pagenumbers#1#2{\def\startpage{#1}\def\finishpage{#2}}
\def\published#1{\def\publishdate{#1}}

\def\received#1{\def\receiveddate{#1}}

\def\accepted#1{\def\accepteddate{#1}}
\def\asciititle#1{\def\theasciititle{#1}}
\def\covertitle#1{\def\thecovertitle{#1}}

\long\def\asciiabstract#1{\long\def\theasciiabstract{#1}}

\let\\\par\let\thelognumber\relax\let\thevolumenumber\relax
\let\thepapernumber\relax\let\thevolumeyear\relax\let\startpage\relax
\let\finishpage\relax\let\publishdate\relax\let\receiveddate\relax
\let\reviseddate\relax\let\accepteddate\relax\let\theasciititle\relax
\let\thecovertitle\relax\let\theasciiauthors\relax
\let\theasciiabstract\relax

\let\theasciiemail\relax

\ifplaintex
\font\logobig=cmssbx10 scaled 3836
\font\logomed=cmssbx10 scaled 2557
\else
\font\logobig=cmssbx10 scaled 4200
\font\logomed=cmssbx10 scaled 2800
\fi

\long\def\makeagttitle{   
\count0=\startpage
\agt\hfill      
\hbox to 45truept{\vbox to 0pt{\vglue -13truept{\logomed A\kern -.37em{\logobig 
T}\kern -.38em G}\vss}\hss}
\break
{\small Volume \thevolumenumber\ (\thevolumeyear)
\startpage--\finishpage\nl
Published: \publishdate}

\vglue .25truein

{\parskip=0pt\leftskip 0pt plus
1fil\def\\{\par\smallskip}{\Large\bf\thetitle}\par\medskip} \vglue
0.05truein

{\parskip=0pt\leftskip 0pt plus 1fil\def\\{\par}{\sc\theauthors}
\par\medskip}%
 
\vglue 0.03truein 

{\small\leftskip 25truept\rightskip 25truept{\bf Abstract}\stdspace\theabstract

{\bf AMS Classification}\stdspace\theprimaryclass
\ifx\thesecondaryclass\relax\else; \thesecondaryclass\fi\par
{\bf Keywords}\stdspace \thekeywords\par}\vglue 7truept

}   

\ifplaintex
\font\phead=cmsl9 scaled 950
\font\pnum=cmbx10 scaled 913
\font\pfoot=cmsl9 scaled 950
\headline{\vbox to 0pt{\vskip -4.5mm\line{\small\phead\ifnum
\count0=\startpage ISSN 1472-2739 (on-line) 1472-2747 (printed)
\hfill {\pnum\folio}\else\ifodd\count0\def\\{ }%
\ifx\theshorttitle\relax\thetitle\else\theshorttitle\fi\hfill{\pnum\folio}
\else\def\\{ and }{\pnum\folio}\hfill\ifx\theshortauthors\relax\theauthors
\else\theshortauthors\fi\fi\fi}\vss}}
\footline{\vbox to 0pt{\vglue 0mm\line{\small\pfoot\ifnum\count0=\startpage
\copyright\ \gtp\hfill\else
\agt, Volume \thevolumenumber\ (\thevolumeyear)\hfill\fi}\vss}}
\else
\font=cmsl9 scaled 1050
\font\lnum=cmbx10 
\font=cmsl9 scaled 1050
\makeatletter
\def\@oddhead{{\small\lhead\ifnum\count0=\startpage ISSN 1472-2739 
(on-line) 1472-2747 (printed)\hfill {\lnum\number\count0}\else\ifodd\count0
\def\\{ }\ifx\theshorttitle\relax \thetitle \else\theshorttitle\fi\hfill
{\lnum\number\count0}\else\def\\{ and }{\lnum\number\count0}
\hfill\ifx\theshortauthors\relax 
\theauthors\else\theshortauthors\fi\fi\fi}}\def\@evenhead{\@oddhead}
\def\@oddfoot{\small\lfoot\ifnum\count0=\startpage\copyright\ \gtp\hfill\else
\agt, Volume \thevolumenumber\ (\thevolumeyear)\hfill\fi}
\def\@evenfoot{\@oddfoot}
\makeatother
\fi
\let\maketitlepage\makeagttitle
\let\makeshorttitle\maketitlepage
\let\maketitle\maketitlepage

\newwrite\gtoutfile
\long\gdef\makeheadfile{  
{\def\\{, }\def\s{ }
\immediate\openout\gtoutfile head.xxx
\immediate\write\gtoutfile{To: math@arxiv.org}
\immediate\write\gtoutfile{Subject: put OR rep NNNNN:ppppp}
\immediate\write\gtoutfile{--text follows this line--}
\immediate\write\gtoutfile{Proxy-for: \ifx\theasciiauthors\relax
\theauthors\else\theasciiauthors\fi\s<\ifx\theasciiemail\relax\theemail\else\theasciiemail\fi>}
\immediate\write\gtoutfile{\noexpand\\}
\immediate\write\gtoutfile{Authors: \ifx\theasciiauthors\relax
\theauthors\else\theasciiauthors\fi}
{\def\\{ }\immediate\write\gtoutfile{Title: \ifx\theasciititle\relax
\thetitle\else\theasciititle\fi}}
\immediate\write\gtoutfile{Subj-class: GT or SG, GR etc}
\immediate\write\gtoutfile{MSC-class: \theprimaryclass\ifx\thesecondaryclass\relax\else, \thesecondaryclass\fi}
\immediate\write\gtoutfile{Journal-ref: Algebr. Geom. Topol. \thevolumenumber\s
(\thevolumeyear) \startpage-\finishpage}
\immediate\write\gtoutfile{Comments: Published by Algebraic and
Geometric Topology at}
\immediate\write\gtoutfile{\s\s\s  http://www.maths.warwick.ac.uk/agt/AGTVol\thevolumenumber/agt-\thevolumenumber-\thepapernumber.abs.html}
\immediate\write\gtoutfile{\noexpand\\}
\immediate\write\gtoutfile{}
\ifx\theasciiabstract\relax
\immediate\write\gtoutfile{\theabstract}\else
\immediate\write\gtoutfile{\theasciiabstract}\fi
\immediate\write\gtoutfile{}
\immediate\write\gtoutfile{\noexpand\\}
\immediate\write\gtoutfile{}
\immediate\closeout\gtoutfile}}  

\def\maketitlepage{\makeagttitle\makeheadfile}
\let\makeshorttitle\maketitlepage
\let\maketitle\maketitlepage

\def\ifplaintex{\expandafter\ifx\csname documentclass\endcsname\relax}

\def\gtp{{\mathsurround=0pt\it $\cal G\mskip-2mu$eometry \&\ 
$\cal T\!\!$opology $\cal P\!$ublications}}  

\def\Addressesr{\bigskip
{\small \parskip 0pt \leftskip 0pt \rightskip 0pt plus 1fil \def\\{\par}
\sl\theaddress\par
\medskip
\rm Email:\stdspace\tt\theemail\hfill\rm Received:\qua\receiveddate \par}}

\def\recd{{\small Received:\qua\receiveddate\ifx\reviseddate\relax
\else\qquad Revised:\qua\reviseddate\fi\par}} 

\def\lognumber#1{\def\thelognumber{#1}}
\def\volumenumber#1{\def\thevolumenumber{#1}}
\def\volumeyear#1{\def\thevolumeyear{#1}}
\def\papernumber#1{\def\thepapernumber{#1}}
\def\pagenumbers#1#2{\def\startpage{#1}\def\finishpage{#2}}
\def\published#1{\def\publishdate{#1}}

\def\received#1{\def\receiveddate{#1}}

\def\accepted#1{\def\accepteddate{#1}}
\def\asciititle#1{\def\theasciititle{#1}}
\def\covertitle#1{\def\thecovertitle{#1}}

\long\def\asciiabstract#1{\long\def\theasciiabstract{#1}}

\let\\\par\let\thelognumber\relax\let\thevolumenumber\relax
\let\thepapernumber\relax\let\thevolumeyear\relax\let\startpage\relax
\let\finishpage\relax\let\publishdate\relax\let\receiveddate\relax
\let\reviseddate\relax\let\accepteddate\relax\let\theasciititle\relax
\let\thecovertitle\relax\let\theasciiauthors\relax
\let\theasciiabstract\relax

\let\theasciiemail\relax

\ifplaintex
\font\logobig=cmssbx10 scaled 3836
\font\logomed=cmssbx10 scaled 2557
\else
\font\logobig=cmssbx10 scaled 4200
\font\logomed=cmssbx10 scaled 2800
\fi

\long\def\makeagttitle{   
\count0=\startpage
\agt\hfill      
\hbox to 45truept{\vbox to 0pt{\vglue -13truept{\logomed A\kern -.37em{\logobig 
T}\kern -.38em G}\vss}\hss}
\break
{\small Volume \thevolumenumber\ (\thevolumeyear)
\startpage--\finishpage\nl
Published: \publishdate}

\vglue .25truein

{\parskip=0pt\leftskip 0pt plus
1fil\def\\{\par\smallskip}{\Large\bf\thetitle}\par\medskip} \vglue
0.05truein

{\parskip=0pt\leftskip 0pt plus 1fil\def\\{\par}{\sc\theauthors}
\par\medskip}%
 
\vglue 0.03truein 

{\small\leftskip 25truept\rightskip 25truept{\bf Abstract}\stdspace\theabstract

{\bf AMS Classification}\stdspace\theprimaryclass
\ifx\thesecondaryclass\relax\else; \thesecondaryclass\fi\par
{\bf Keywords}\stdspace \thekeywords\par}\vglue 7truept

}   

\ifplaintex
\font\phead=cmsl9 scaled 950
\font\pnum=cmbx10 scaled 913
\font\pfoot=cmsl9 scaled 950
\headline{\vbox to 0pt{\vskip -4.5mm\line{\small\phead\ifnum
\count0=\startpage ISSN 1472-2739 (on-line) 1472-2747 (printed)
\hfill {\pnum\folio}\else\ifodd\count0\def\\{ }%
\ifx\theshorttitle\relax\thetitle\else\theshorttitle\fi\hfill{\pnum\folio}
\else\def\\{ and }{\pnum\folio}\hfill\ifx\theshortauthors\relax\theauthors
\else\theshortauthors\fi\fi\fi}\vss}}
\footline{\vbox to 0pt{\vglue 0mm\line{\small\pfoot\ifnum\count0=\startpage
\copyright\ \gtp\hfill\else
\agt, Volume \thevolumenumber\ (\thevolumeyear)\hfill\fi}\vss}}
\else
\font=cmsl9 scaled 1050
\font\lnum=cmbx10 
\font=cmsl9 scaled 1050
\makeatletter
\def\@oddhead{{\small\lhead\ifnum\count0=\startpage ISSN 1472-2739 
(on-line) 1472-2747 (printed)\hfill {\lnum\number\count0}\else\ifodd\count0
\def\\{ }\ifx\theshorttitle\relax \thetitle \else\theshorttitle\fi\hfill
{\lnum\number\count0}\else\def\\{ and }{\lnum\number\count0}
\hfill\ifx\theshortauthors\relax 
\theauthors\else\theshortauthors\fi\fi\fi}}\def\@evenhead{\@oddhead}
\def\@oddfoot{\small\lfoot\ifnum\count0=\startpage\copyright\ \gtp\hfill\else
\agt, Volume \thevolumenumber\ (\thevolumeyear)\hfill\fi}
\def\@evenfoot{\@oddfoot}
\makeatother
\fi
\let\maketitlepage\makeagttitle
\let\makeshorttitle\maketitlepage
\let\maketitle\maketitlepage

\newwrite\gtoutfile
\long\gdef\makeheadfile{  
{\def\\{, }\def\s{ }
\immediate\openout\gtoutfile head.xxx
\immediate\write\gtoutfile{To: math@arxiv.org}
\immediate\write\gtoutfile{Subject: put OR rep NNNNN:ppppp}
\immediate\write\gtoutfile{--text follows this line--}
\immediate\write\gtoutfile{Proxy-for: \ifx\theasciiauthors\relax
\theauthors\else\theasciiauthors\fi\s<\ifx\theasciiemail\relax\theemail\else\theasciiemail\fi>}
\immediate\write\gtoutfile{\noexpand\\}
\immediate\write\gtoutfile{Authors: \ifx\theasciiauthors\relax
\theauthors\else\theasciiauthors\fi}
{\def\\{ }\immediate\write\gtoutfile{Title: \ifx\theasciititle\relax
\thetitle\else\theasciititle\fi}}
\immediate\write\gtoutfile{Subj-class: GT or SG, GR etc}
\immediate\write\gtoutfile{MSC-class: \theprimaryclass\ifx\thesecondaryclass\relax\else, \thesecondaryclass\fi}
\immediate\write\gtoutfile{Journal-ref: Algebr. Geom. Topol. \thevolumenumber\s
(\thevolumeyear) \startpage-\finishpage}
\immediate\write\gtoutfile{Comments: Published by Algebraic and
Geometric Topology at}
\immediate\write\gtoutfile{\s\s\s  http://www.maths.warwick.ac.uk/agt/AGTVol\thevolumenumber/agt-\thevolumenumber-\thepapernumber.abs.html}
\immediate\write\gtoutfile{\noexpand\\}
\immediate\write\gtoutfile{}
\ifx\theasciiabstract\relax
\immediate\write\gtoutfile{\theabstract}\else
\immediate\write\gtoutfile{\theasciiabstract}\fi
\immediate\write\gtoutfile{}
\immediate\write\gtoutfile{\noexpand\\}
\immediate\write\gtoutfile{}
\immediate\closeout\gtoutfile}}  

\def\maketitlepage{\makeagttitle\makeheadfile}
\let\makeshorttitle\maketitlepage
\let\maketitle\maketitlepage

\def\ifplaintex{\expandafter\ifx\csname documentclass\endcsname\relax}

\def\gtp{{\mathsurround=0pt\it $\cal G\mskip-2mu$eometry \&\ 
$\cal T\!\!$opology $\cal P\!$ublications}}  

\def\Addressesr{\bigskip
{\small \parskip 0pt \leftskip 0pt \rightskip 0pt plus 1fil \def\\{\par}
\sl\theaddress\par
\medskip
\rm Email:\stdspace\tt\theemail\hfill\rm Received:\qua\receiveddate \par}}

\def\recd{{\small Received:\qua\receiveddate\ifx\reviseddate\relax
\else\qquad Revised:\qua\reviseddate\fi\par}} 

\def\lognumber#1{\def\thelognumber{#1}}
\def\volumenumber#1{\def\thevolumenumber{#1}}
\def\volumeyear#1{\def\thevolumeyear{#1}}
\def\papernumber#1{\def\thepapernumber{#1}}
\def\pagenumbers#1#2{\def\startpage{#1}\def\finishpage{#2}}
\def\published#1{\def\publishdate{#1}}

\def\received#1{\def\receiveddate{#1}}

\def\accepted#1{\def\accepteddate{#1}}
\def\asciititle#1{\def\theasciititle{#1}}
\def\covertitle#1{\def\thecovertitle{#1}}

\long\def\asciiabstract#1{\long\def\theasciiabstract{#1}}

\let\\\par\let\thelognumber\relax\let\thevolumenumber\relax
\let\thepapernumber\relax\let\thevolumeyear\relax\let\startpage\relax
\let\finishpage\relax\let\publishdate\relax\let\receiveddate\relax
\let\reviseddate\relax\let\accepteddate\relax\let\theasciititle\relax
\let\thecovertitle\relax\let\theasciiauthors\relax
\let\theasciiabstract\relax

\let\theasciiemail\relax

\ifplaintex
\font\logobig=cmssbx10 scaled 3836
\font\logomed=cmssbx10 scaled 2557
\else
\font\logobig=cmssbx10 scaled 4200
\font\logomed=cmssbx10 scaled 2800
\fi

\long\def\makeagttitle{   
\count0=\startpage
\agt\hfill      
\hbox to 45truept{\vbox to 0pt{\vglue -13truept{\logomed A\kern -.37em{\logobig 
T}\kern -.38em G}\vss}\hss}
\break
{\small Volume \thevolumenumber\ (\thevolumeyear)
\startpage--\finishpage\nl
Published: \publishdate}

\vglue .25truein

{\parskip=0pt\leftskip 0pt plus
1fil\def\\{\par\smallskip}{\Large\bf\thetitle}\par\medskip} \vglue
0.05truein

{\parskip=0pt\leftskip 0pt plus 1fil\def\\{\par}{\sc\theauthors}
\par\medskip}%
 
\vglue 0.03truein 

{\small\leftskip 25truept\rightskip 25truept{\bf Abstract}\stdspace\theabstract

{\bf AMS Classification}\stdspace\theprimaryclass
\ifx\thesecondaryclass\relax\else; \thesecondaryclass\fi\par
{\bf Keywords}\stdspace \thekeywords\par}\vglue 7truept

}   

\ifplaintex
\font\phead=cmsl9 scaled 950
\font\pnum=cmbx10 scaled 913
\font\pfoot=cmsl9 scaled 950
\headline{\vbox to 0pt{\vskip -4.5mm\line{\small\phead\ifnum
\count0=\startpage ISSN 1472-2739 (on-line) 1472-2747 (printed)
\hfill {\pnum\folio}\else\ifodd\count0\def\\{ }%
\ifx\theshorttitle\relax\thetitle\else\theshorttitle\fi\hfill{\pnum\folio}
\else\def\\{ and }{\pnum\folio}\hfill\ifx\theshortauthors\relax\theauthors
\else\theshortauthors\fi\fi\fi}\vss}}
\footline{\vbox to 0pt{\vglue 0mm\line{\small\pfoot\ifnum\count0=\startpage
\copyright\ \gtp\hfill\else
\agt, Volume \thevolumenumber\ (\thevolumeyear)\hfill\fi}\vss}}
\else
\font=cmsl9 scaled 1050
\font\lnum=cmbx10 
\font=cmsl9 scaled 1050
\makeatletter
\def\@oddhead{{\small\lhead\ifnum\count0=\startpage ISSN 1472-2739 
(on-line) 1472-2747 (printed)\hfill {\lnum\number\count0}\else\ifodd\count0
\def\\{ }\ifx\theshorttitle\relax \thetitle \else\theshorttitle\fi\hfill
{\lnum\number\count0}\else\def\\{ and }{\lnum\number\count0}
\hfill\ifx\theshortauthors\relax 
\theauthors\else\theshortauthors\fi\fi\fi}}\def\@evenhead{\@oddhead}
\def\@oddfoot{\small\lfoot\ifnum\count0=\startpage\copyright\ \gtp\hfill\else
\agt, Volume \thevolumenumber\ (\thevolumeyear)\hfill\fi}
\def\@evenfoot{\@oddfoot}
\makeatother
\fi
\let\maketitlepage\makeagttitle
\let\makeshorttitle\maketitlepage
\let\maketitle\maketitlepage

\newwrite\gtoutfile
\long\gdef\makeheadfile{  
{\def\\{, }\def\s{ }
\immediate\openout\gtoutfile head.xxx
\immediate\write\gtoutfile{To: math@arxiv.org}
\immediate\write\gtoutfile{Subject: put OR rep NNNNN:ppppp}
\immediate\write\gtoutfile{--text follows this line--}
\immediate\write\gtoutfile{Proxy-for: \ifx\theasciiauthors\relax
\theauthors\else\theasciiauthors\fi\s<\ifx\theasciiemail\relax\theemail\else\theasciiemail\fi>}
\immediate\write\gtoutfile{\noexpand\\}
\immediate\write\gtoutfile{Authors: \ifx\theasciiauthors\relax
\theauthors\else\theasciiauthors\fi}
{\def\\{ }\immediate\write\gtoutfile{Title: \ifx\theasciititle\relax
\thetitle\else\theasciititle\fi}}
\immediate\write\gtoutfile{Subj-class: GT or SG, GR etc}
\immediate\write\gtoutfile{MSC-class: \theprimaryclass\ifx\thesecondaryclass\relax\else, \thesecondaryclass\fi}
\immediate\write\gtoutfile{Journal-ref: Algebr. Geom. Topol. \thevolumenumber\s
(\thevolumeyear) \startpage-\finishpage}
\immediate\write\gtoutfile{Comments: Published by Algebraic and
Geometric Topology at}
\immediate\write\gtoutfile{\s\s\s  http://www.maths.warwick.ac.uk/agt/AGTVol\thevolumenumber/agt-\thevolumenumber-\thepapernumber.abs.html}
\immediate\write\gtoutfile{\noexpand\\}
\immediate\write\gtoutfile{}
\ifx\theasciiabstract\relax
\immediate\write\gtoutfile{\theabstract}\else
\immediate\write\gtoutfile{\theasciiabstract}\fi
\immediate\write\gtoutfile{}
\immediate\write\gtoutfile{\noexpand\\}
\immediate\write\gtoutfile{}
\immediate\closeout\gtoutfile}}  

\def\maketitlepage{\makeagttitle\makeheadfile}
\let\makeshorttitle\maketitlepage
\let\maketitle\maketitlepage

\def\ifplaintex{\expandafter\ifx\csname documentclass\endcsname\relax}

\def\gtp{{\mathsurround=0pt\it $\cal G\mskip-2mu$eometry \&\ 
$\cal T\!\!$opology $\cal P\!$ublications}}  

\def\Addressesr{\bigskip
{\small \parskip 0pt \leftskip 0pt \rightskip 0pt plus 1fil \def\\{\par}
\sl\theaddress\par
\medskip
\rm Email:\stdspace\tt\theemail\hfill\rm Received:\qua\receiveddate \par}}

\def\recd{{\small Received:\qua\receiveddate\ifx\reviseddate\relax
\else\qquad Revised:\qua\reviseddate\fi\par}} 

\def\lognumber#1{\def\thelognumber{#1}}
\def\volumenumber#1{\def\thevolumenumber{#1}}
\def\volumeyear#1{\def\thevolumeyear{#1}}
\def\papernumber#1{\def\thepapernumber{#1}}
\def\pagenumbers#1#2{\def\startpage{#1}\def\finishpage{#2}}
\def\published#1{\def\publishdate{#1}}

\def\received#1{\def\receiveddate{#1}}

\def\accepted#1{\def\accepteddate{#1}}
\def\asciititle#1{\def\theasciititle{#1}}
\def\covertitle#1{\def\thecovertitle{#1}}

\long\def\asciiabstract#1{\long\def\theasciiabstract{#1}}

\let\\\par\let\thelognumber\relax\let\thevolumenumber\relax
\let\thepapernumber\relax\let\thevolumeyear\relax\let\startpage\relax
\let\finishpage\relax\let\publishdate\relax\let\receiveddate\relax
\let\reviseddate\relax\let\accepteddate\relax\let\theasciititle\relax
\let\thecovertitle\relax\let\theasciiauthors\relax
\let\theasciiabstract\relax

\let\theasciiemail\relax

\ifplaintex
\font\logobig=cmssbx10 scaled 3836
\font\logomed=cmssbx10 scaled 2557
\else
\font\logobig=cmssbx10 scaled 4200
\font\logomed=cmssbx10 scaled 2800
\fi

\long\def\makeagttitle{   
\count0=\startpage
\agt\hfill      
\hbox to 45truept{\vbox to 0pt{\vglue -13truept{\logomed A\kern -.37em{\logobig 
T}\kern -.38em G}\vss}\hss}
\break
{\small Volume \thevolumenumber\ (\thevolumeyear)
\startpage--\finishpage\nl
Published: \publishdate}

\vglue .25truein

{\parskip=0pt\leftskip 0pt plus
1fil\def\\{\par\smallskip}{\Large\bf\thetitle}\par\medskip} \vglue
0.05truein

{\parskip=0pt\leftskip 0pt plus 1fil\def\\{\par}{\sc\theauthors}
\par\medskip}%
 
\vglue 0.03truein 

{\small\leftskip 25truept\rightskip 25truept{\bf Abstract}\stdspace\theabstract

{\bf AMS Classification}\stdspace\theprimaryclass
\ifx\thesecondaryclass\relax\else; \thesecondaryclass\fi\par
{\bf Keywords}\stdspace \thekeywords\par}\vglue 7truept

}   

\ifplaintex
\font\phead=cmsl9 scaled 950
\font\pnum=cmbx10 scaled 913
\font\pfoot=cmsl9 scaled 950
\headline{\vbox to 0pt{\vskip -4.5mm\line{\small\phead\ifnum
\count0=\startpage ISSN 1472-2739 (on-line) 1472-2747 (printed)
\hfill {\pnum\folio}\else\ifodd\count0\def\\{ }%
\ifx\theshorttitle\relax\thetitle\else\theshorttitle\fi\hfill{\pnum\folio}
\else\def\\{ and }{\pnum\folio}\hfill\ifx\theshortauthors\relax\theauthors
\else\theshortauthors\fi\fi\fi}\vss}}
\footline{\vbox to 0pt{\vglue 0mm\line{\small\pfoot\ifnum\count0=\startpage
\copyright\ \gtp\hfill\else
\agt, Volume \thevolumenumber\ (\thevolumeyear)\hfill\fi}\vss}}
\else
\font=cmsl9 scaled 1050
\font\lnum=cmbx10 
\font=cmsl9 scaled 1050
\makeatletter
\def\@oddhead{{\small\lhead\ifnum\count0=\startpage ISSN 1472-2739 
(on-line) 1472-2747 (printed)\hfill {\lnum\number\count0}\else\ifodd\count0
\def\\{ }\ifx\theshorttitle\relax \thetitle \else\theshorttitle\fi\hfill
{\lnum\number\count0}\else\def\\{ and }{\lnum\number\count0}
\hfill\ifx\theshortauthors\relax 
\theauthors\else\theshortauthors\fi\fi\fi}}\def\@evenhead{\@oddhead}
\def\@oddfoot{\small\lfoot\ifnum\count0=\startpage\copyright\ \gtp\hfill\else
\agt, Volume \thevolumenumber\ (\thevolumeyear)\hfill\fi}
\def\@evenfoot{\@oddfoot}
\makeatother
\fi
\let\maketitlepage\makeagttitle
\let\makeshorttitle\maketitlepage
\let\maketitle\maketitlepage

\newwrite\gtoutfile
\long\gdef\makeheadfile{  
{\def\\{, }\def\s{ }
\immediate\openout\gtoutfile head.xxx
\immediate\write\gtoutfile{To: math@arxiv.org}
\immediate\write\gtoutfile{Subject: put OR rep NNNNN:ppppp}
\immediate\write\gtoutfile{--text follows this line--}
\immediate\write\gtoutfile{Proxy-for: \ifx\theasciiauthors\relax
\theauthors\else\theasciiauthors\fi\s<\ifx\theasciiemail\relax\theemail\else\theasciiemail\fi>}
\immediate\write\gtoutfile{\noexpand\\}
\immediate\write\gtoutfile{Authors: \ifx\theasciiauthors\relax
\theauthors\else\theasciiauthors\fi}
{\def\\{ }\immediate\write\gtoutfile{Title: \ifx\theasciititle\relax
\thetitle\else\theasciititle\fi}}
\immediate\write\gtoutfile{Subj-class: GT or SG, GR etc}
\immediate\write\gtoutfile{MSC-class: \theprimaryclass\ifx\thesecondaryclass\relax\else, \thesecondaryclass\fi}
\immediate\write\gtoutfile{Journal-ref: Algebr. Geom. Topol. \thevolumenumber\s
(\thevolumeyear) \startpage-\finishpage}
\immediate\write\gtoutfile{Comments: Published by Algebraic and
Geometric Topology at}
\immediate\write\gtoutfile{\s\s\s  http://www.maths.warwick.ac.uk/agt/AGTVol\thevolumenumber/agt-\thevolumenumber-\thepapernumber.abs.html}
\immediate\write\gtoutfile{\noexpand\\}
\immediate\write\gtoutfile{}
\ifx\theasciiabstract\relax
\immediate\write\gtoutfile{\theabstract}\else
\immediate\write\gtoutfile{\theasciiabstract}\fi
\immediate\write\gtoutfile{}
\immediate\write\gtoutfile{\noexpand\\}
\immediate\write\gtoutfile{}
\immediate\closeout\gtoutfile}}  

\def\maketitlepage{\makeagttitle\makeheadfile}
\let\makeshorttitle\maketitlepage
\let\maketitle\maketitlepage

\lognumber{21}
\volumenumber{2}
\volumeyear{2002}
\papernumber{21}
\published{29 May 2002}
\pagenumbers{433}{447}
\received{25 March 2002}
\accepted{28 May 2002}

\usepackage{epsf, amsmath, amsthm, amssymb}

\input epsf
\def\relabelbox{%
  \hbox\bgroup%
}%
\def\endrelabelbox{%
\egroup%
}%
\def\relabel #1#2 {%
  \special{ps:/a {} def}%
  \smash{\rlap{#2}}%
}%
\def\adjustrelabel <#1,#2> #3#4 {%
  \special{ps:/a {} def}%
  \smash{\rlap{\kern #1 \raise #2\hbox{#4}}}%
}%
\def\extralabel <#1,#2> #3 {\smash{\rlap{\kern #1 \raise #2\hbox{#3}}}}%

\let\mathscr\mathcal

\begin{document}

\newtheorem{thm}{Theorem}[section]
\newtheorem{lem}[thm]{Lemma}
\newtheorem{cor}[thm]{Corollary}
\newtheorem{conj}[thm]{Conjecture}
\newtheorem{qn}[thm]{Question}

\theoremstyle{definition}
\newtheorem{defn}[thm]{Definition}

\theoremstyle{remark}
\newtheorem*{rmk}{Remark}
\newtheorem*{exa}{Example}

\def\R{\mathbb R}
\def\Z{\mathbb Z}
\def\CP{\mathbb {CP}}
\def\H{\mathscr H}
\def\h2{{\mathbb H}^2}
\def\F{\mathscr F}
\def\E{\mathscr E}
\def\M{\mathscr M}
\def\C{\mathscr C}
\def\Q{\mathbb Q}
\def\S{\mathscr S}
\def\T{\mathscr T}
\def\G{\mathscr G}
\def\L{\mathcal L}
\def\O{\mathscr O}
\def\BG{{{\rm B}\Gamma}}
\def\BO{{{\rm BO}}}
\def\homeo{\text{Homeo}}
\def\u{{\text{univ}}}
\def\f{{\text{fix}}}
\def\inte{{\text{int}}}
\def\id{{\text{id}}}
\def\stb{{\text{stab}}}
\def\sli{{\text{sl}}}
\def\til{\widetilde}
\def\I{\mathscr I}
\def\N{\mathbb N}

\title{Every orientable $3$--manifold is a $\BG$}
\covertitle{Every orientable $3$--manifold is a ${\noexpand\rm B}\Gamma$}
\asciititle{Every orientable 3-manifold is a B\Gamma}

\author{Danny Calegari}

\address{Department of Mathematics, Harvard University\\Cambridge MA, 02138, USA}
\email{dannyc@math.harvard.edu}

\begin{abstract}
We show that every orientable $3$--manifold is a classifying space
$\BG$ where $\Gamma$ is a groupoid of germs of homeomorphisms of
$\R$. This follows by showing that every orientable $3$--manifold $M$
admits a codimension one foliation $\F$ such that the holonomy cover
of every leaf is contractible. The $\F$ we construct can be taken to
be $C^1$ but not $C^2$.  The existence of such an $\F$ answers
positively a question posed by Tsuboi (see \cite{tT01}), but leaves
open the question of whether $M = \BG$ for some $C^\infty$ groupoid
$\Gamma$.
\end{abstract}

\asciiabstract{We show that every orientable 3-manifold is a
classifying space B\Gamma where \Gamma is a groupoid of germs of
homeomorphisms of R. This follows by showing that every orientable
3-manifold M admits a codimension one foliation F such that the
holonomy cover of every leaf is contractible. The F we construct can
be taken to be C^1 but not C^2.  The existence of such an F answers
positively a question posed by Tsuboi [Classifying spaces for groupoid
structures, notes from minicourse at PUC, Rio de Janeiro (2001)], but
leaves open the question of whether M = B\Gamma for some C^\infty
groupoid \Gamma.}

\maketitle

\section{Acknowledgements}
In August 2001, the P.U.C. in Rio de Janeiro held a conference on Foliations and
Dynamics aimed at bringing together traditional foliators and $3$--manifold
topologists. At this conference, Takashi Tsuboi posed the question of the existence
of typical foliations on $3$--manifolds. I would very much like to thank Takashi
for posing this question, for reading an early draft of this paper and
catching numerous errors, and for introducting me to the beautiful
subject of classical foliation theory. I'd also like to thank Curt McMullen for a
useful conversation which helped clarify some analytic questions for me.

\section{Typical foliations}

Let $M$ be a smooth $n$--manifold, and $\F$ a smooth codimension $p$ foliation of $M$.
The tangent bundle $T\F$ is a subbundle of $TM$, with orthogonal complement the
normal bundle $\nu\F$. Choose some small $\epsilon > 0$, and let
$E_\epsilon \subset \nu\F$ be the $B^p$ bundle over $M$ consisting of vectors of $\nu\F$
of length $\le \epsilon$. For sufficiently small $\epsilon$, for each
leaf $\lambda$ of $\F$, the restriction of $E_\epsilon$ 
to $\lambda$ maps via the exponential map to a (typically immersed)
tubular neighborhood of $\lambda$. The distribution
$T\F$ can be pulled back via this exponential map to a (generally incomplete) connection
on $E_\epsilon$. For each point $p \in \lambda$, the ball $E_\epsilon(p)$ parameterizes the
leaves of $\F$ near $p$. If $\gamma:I \to \lambda$ is a path,
the holonomy of the connection on $E_\epsilon$ defines a diffeomorphism from a neighborhood of
$0$ in $E_\epsilon(\gamma(0))$ to a neighborhood of $0$ in $E_\epsilon(\gamma(1))$.
The germ of this map is called the {\em holonomy of $\F$ along $\gamma$},
and depends only on the homotopy class of $\gamma$ rel.
endpoints in $\lambda$. In particular, for any
$p \in \lambda$, there is a homomorphism from $\pi_1(\lambda,p)$ to the group of
germs at $0$ of diffeomorphisms of $E_\epsilon(p)$ to itself which fix $0$. The {\em holonomy
cover} of $\lambda$ is the covering space $\widetilde{\lambda}$ of $\lambda$ corresponding
to the kernel of this homomorphism.

If $M$ and $\F$ are not necessarily smooth, there is still a well--defined holonomy
homomorphism from $\pi_1(\lambda,p)$ to the group of germs at $p$ of homeomorphisms of
a transversal $\tau$ to $\F$ through $p$ to itself fixing $p$. Therefore we can define
the holonomy cover of a leaf in this case too.

\begin{defn}
A foliation $\F$ is {\em typical} if for each leaf $\lambda$, the
holonomy cover $\widetilde{\lambda}$ of $\lambda$ is contractible.
\end{defn}

If we fix a finite collection $\lbrace \tau_i \rbrace$
of transversals to $\F$, let $\pi_1(\F,\lbrace \tau_i\rbrace)$
denote the groupoid of homotopy classes of paths rel. endpoints
tangent to $\F$, which start and end on some pair of transversals
$\tau_i,\tau_j \in \lbrace \tau_i \rbrace$. This groupoid
maps by the holonomy homomorphism
to the groupoid of germs of homeomorphisms between subsets of the $\tau_i$.
If the union of the $\tau_i$ are a {\em total transversal} --- i.e.\ if they
intersect every leaf of $\F$ ---
we call the image of $\pi_1(\F,\lbrace \tau_i \rbrace)$ the
{\em holonomy groupoid of the foliation on $\lbrace \tau_i \rbrace$}. By
embedding the $\tau_i$ disjointly in $\R^p$, we can think of the holonomy
groupoid as a groupoid of germs of homeomorphisms of $\R^p$.

We are especially interested, throughout most of this article, with surface
foliations of $3$--manifolds. A surface foliation without sphere or projective plane
leaves is typical iff the fundamental group $\pi_1(\lambda)$ of each leaf injects into the
holonomy groupoid of $\F$.
For a surface foliation of a $3$--manifold, the holonomy groupoid $\Gamma$
is a groupoid of germs of homeomorphisms of $\R$.
If $\F$ is $C^\infty$, the holonomy groupoid $\Gamma$ is a groupoid of
germs of diffeomorphisms of $\R$. For a more leisurely discussion of
holonomy groupoids, see \cite{pM88}.

The significance of typical foliations is explained by the following
theorem of Haefliger:

\begin{thm}[Haefliger]\label{typical_classifies}
Let $\F$ be a typical foliation of $M$. Then $M$ is homotopic to $\BG$
where $\Gamma$ is the holonomy groupoid of $\F$ (on any total transversal).
\end{thm}

Here $\BG$ is the {\em classifying space} of the groupoid $\Gamma$.
For a construction of the classifying space of a groupoid, see \cite{gS78}. 
Note that for $\Gamma$ a groupoid of germs of homeomorphisms of $\R^p$,
$\Gamma$ is a {\em topological groupoid}, although the topology is not necessarily what
one might expect. For $\phi \in \homeo(\R^p)$ with $\phi(x) = y$, we denote the
germ of $\phi$ at $x$ by $\phi_x$.
We give $\homeo(\R^p)$ the {\em discrete} topology and $\R^p$ the
usual topology, and think of $\Gamma$
as a quotient of the subspace
$$\lbrace (\phi,x,y) \in \homeo(\R^p) \times \R^p \times \R^p: \phi_x \in \Gamma,
\phi(x) = y \rbrace$$
by the equivalence relation:
$$(\phi,x,y) \sim (\psi,x,y) \text{ if } \phi_x = \psi_x$$
For $X$ a topological space, the set of homotopy classes of
Haefliger structures on $X$
with coefficients in a groupoid $\Gamma$ are classified by homotopy
classes of maps $[X,\BG]$. If $\Gamma$ is the groupoid of germs of $C^r$ diffeomorphisms
of $\R^p$ for some $r \in \lbrace 1,\dots,\infty \rbrace$, there is a natural ``forgetful''
map
$$\pi:\BG \to \BO(p,\R)$$
which remembers only the linear part of each $\gamma \in \Gamma$.

Let $\xi$ be the tautological bundle over $\BO(p,\R)$.
If $M$ is a smooth manifold, Thurston showed (\cite{wT74b},\cite{wT76}) that for every
$\alpha:M \to \BG$ for which $(\pi\circ\alpha)^*\xi$ is isomorphic to
a sub-bundle $\zeta \subset TM$,
there is a $C^r$ foliation $\F$ of $M$ with $\nu\F$ homotopic to $\zeta$,
and with Haefliger structure homotopic to $\alpha^*\H$, for $\H$ the tautological Haefliger
structure on $\BG$. See \cite{tT01} and \cite{tT90}
for the definition of a Haefliger structure and for
more details.

\begin{exa}
If $\F$ is a $p$--dimensional foliation with every leaf homeomorphic to
$\R^p$, then $\F$ is typical. Some familiar examples include:
\begin{enumerate}
\item{Kronecker foliations. Let $\pi$ be a $p$--dimensional subspace
  of $\R^n$ which intersects the integer lattice $\Z^n$ only at the origin.
Let $\til{\F}$ be the foliation of $\R^n$ by planes parallel to $\pi$.
Then $\til{\F}$ descends to a foliation $\F$ of $T^n = \R^n/\Z^n$,
all of whose leaves are $p$--planes.}
\item{Hilbert modular surfaces. Let $\O$ denote the ring of integers in
a real quadratic extension $K$ of $\Q$, and let $e_1,e_2:K \to \R$ denote the two
real embeddings. Let $\Gamma$ denote the group $PSL(2,\O)$, and let $\widehat{\Gamma}$
denote a finite index torsion--free subgroup. Then $e_i$ induces an embedding
$$\widehat{e_i}:\widehat{\Gamma} \to PSL(2,\R)$$\eject
for $i=1,2$, and we can let $M$ be the quotient of $\h2 \times \h2$ by 
$\widehat{\Gamma}$, where
$$\gamma(x,y) = (\widehat{e_1}(\gamma)(x),\widehat{e_2}(\gamma)(y))$$
The foliation of $\h2 \times \h2$ by hyperbolic planes $\h2 \times \text{point}$
descends to a foliation of $M$ by hyperbolic planes.}
\item{The Seifert conjecture. Let $M$ be a $3$--manifold, and $X$ a nowhere zero
vector field on $M$. Then Schweitzer (\cite{pS74})
showed how to modify $X$ to $X'$ so that it possesses no closed
integral curves; i.e.\ the integral foliation of $X'$ has all leaves homeomorphic to
$\R$.}
\end{enumerate}
\end{exa}

Many $3$--manifolds are already $\BG$ for $\Gamma$ a {\em group} of germs of
homeomorphisms of $\R$ at some global fixed point.

\begin{exa}
The group $\homeo^+(\R)$ of orientation preserving homeomorphisms of $\R$
is isomorphic by a suspension trick to a subgroup of the
group of germs at $0$ of elements of $\homeo^+(\R)$ which fix $0$. A finitely generated group
$\Gamma$ is isomorphic to a subgroup of $\homeo^+(\R)$ iff $\Gamma$ is {\em left orderable}.
A closed $3$--manifold $M$ is a $K(\pi_1(M),1)$ iff $\pi_1(M)$ is infinite and $\pi_2(M)$ is
trivial; if $\pi_1(M)$ is also left orderable, then $\pi_1(M)$ is isomorphic to a
group of germs of homeomorphisms of $\R$, and $M$ is a classifying space for this group.
\begin{enumerate}
\item If $\Gamma$ is a finitely generated group which is
{\it locally indicable} --- i.e.\ every nontrivial finitely generated
subgroup admits a surjective homomorphism to $\Z$ ---
then $\Gamma$ is isomorphic to a subgroup of $\homeo^+(\R)$.
For $\Gamma$ the fundamental group of an irreducible
$3$--manifold $M$, this happens when $H^1(M)$ is nontrivial.
This observation is due to Boyer, Rolfsen and Wiest \cite{BRW}.

\item{If $M$ is atoroidal and
admits an orientable and co--orientable taut
foliation $\F$, $\pi_1(M)$ admits a faithful representation in
$\homeo^+(S^1)$. If the Euler class of $T\F$ is trivial, this representation lifts to
a subgroup of $\homeo^+(\R)$. If $H^1(M)$ is infinite, $\pi_1(M)$ is locally indicable
as above.  Otherwise, the Euler class of $T\F$ is torsion,
and a finite index subgroup of $\pi_1(M)$ admits a
faithful representation in $\homeo^+(\R)$. See \cite{CD02} for details.}
\end{enumerate}
\end{exa}

Roberts, Shareshian and Stein showed (\cite{RSS}) that there are examples of
hyperbolic $3$--manifold groups which are not left--orderable. More examples and a
detailed analysis are given in \cite{CD02}. 

On the other hand, these examples cannot in general be improved to $C^1$. If
$\Gamma$ is a finitely generated group of germs at $0$ of
$C^1$ homeomorphisms of $\R$ fixing $0$, then $H^1(\Gamma)$ is nontrivial, by
Thurston's stability theorem \cite{wT74a}. As the smoothness increases, our knowledge
shrinks very rapidly. In fact, it is not even known if
$\pi_1(\Sigma_g)$ admits a faithful representation in the group of germs at $0$ of
$C^2$ homeomorphisms of $\R$ fixing $0$, for $\Sigma_g$ the closed orientable
surface of genus $g \ge 2$.

\section{Open book decompositions}

To build our foliations we will appeal to the following simple
structure theorem.

\begin{thm}[Myers]\label{open_book}
Let $M$ be a closed orientable $3$--manifold. Then there is a decomposition
of $M$ as $K \coprod N$ where $K$ is a tubular neighborhood of a knot,
and $N$ is a surface bundle over $S^1$ with fiber the punctured
surface $\Sigma$, such that the boundary curve $\partial \Sigma$ is
not a meridian for $K$. 
\end{thm}

Such a structure is called a {\em open book decomposition with connected
binding}. Such a structure exists for any $3$--manifold by Myers \cite{rM78}.
Note that $\Sigma$ may be chosen to be a Seifert surface for $K$ so
that the boundary curve of $\partial \Sigma$ is not a meridian.

The existence of an open book decomposition for an orientable $3$--manifold
follows from structure theorems of Heegaard and Lickorish (\cite{rL62}). First, $M$ has a
presentation as a union of two handlebodies --- that is, it admits a Heegaard
decomposition. Any two Heegaard decompositions of the same genus
differ by cutting along the splitting surface and regluing by an automorphism
of the surface. Lickorish showed that the automorphism group of the surface
(up to isotopy) is generated by Dehn twists in a canonical set of curves; translating
this into the $3$--dimensional context, and starting from the standard genus $g$
splitting of $S^3$, one sees that any $3$--manifold $M$ of Heegaard genus $g$ is
obtained by integer surgery on a link in $S^3$ of the form illustrated in
figure 1, copied from \cite{dR76}.

\begin{figure}[ht]
\centerline{\relabelbox\small \epsfxsize 3.7truein
\epsfbox{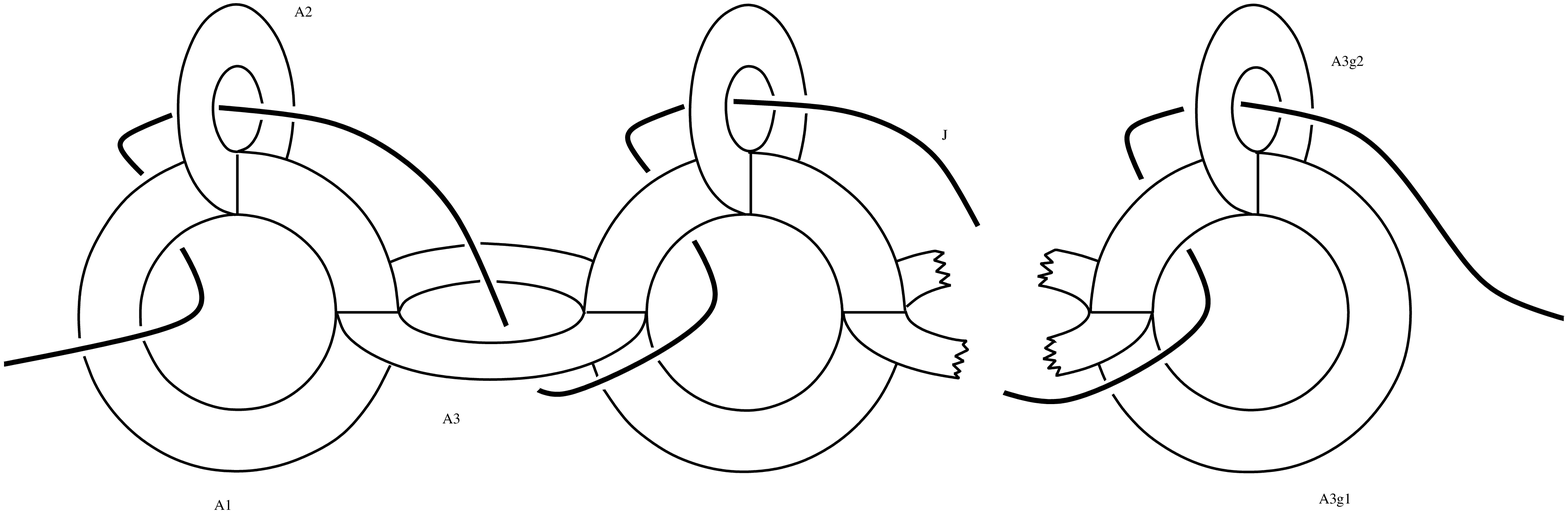}
\relabel {J}{$J$}
\adjustrelabel <0pt,-5pt> {A1}{$A_1$}
\relabel {A2}{$A_2$}
\adjustrelabel <-12pt,-6pt> {A3g1}{$A_{3g-2}$}
\adjustrelabel <0pt,10pt> {A3g2}{$A_{3g-1}$}
\endrelabelbox}
\caption{Any orientable $3$--manifold of Heegaard genus $g$ is obtained
by integral surgery along a collection of components parallel to the cores
of the annuli $A_i$ in the figure. These links form a braid wrapping once around
the curve $J$.}
\end{figure}

$M$ has the structure of an open book decomposition with one component for the
binding for each component of the surgery presentation, and one component coming
from the knot labelled $J$ in figure 1. The fiber is a many--punctured disk.
The number of components of the binding can be reduced at the cost of raising the
genus of the fiber, until there is a single binding component.

In what follows, the fact that our open book decompositions have connected
bindings is not an essential fact, but it simplifies the exposition
somewhat.

Since $N$ is a surface bundle over a circle, it is determined by the
class of the monodromy $\phi:\Sigma \to \Sigma$ up to homotopy. Such
a class contains many smooth representatives.

\section{Foliations of surface bundles}

We inherit the notation $M,N,K,\Sigma$ from the last section.
Let $I$ denote the closed unit interval.

\begin{defn}
A representation $\rho:\pi_1(\Sigma) \to \homeo(I)$ is {\em generic}
if for each nontrivial $\alpha \in \pi_1(\Sigma)$
the fixed points of $\rho(\alpha)$ in the interior of $I$ are
{\em isolated}.
A generic representation whose image consists of $C^\infty$ diffeomorphisms which are
infinitely tangent to the identity at the two endpoints is called {\em well--tapered}.
\end{defn}

The interior of $I$ is homeomorphic to $\R$, and every self--homeomorphism of 
$\inte(I)$ extends uniquely to a self--homeomorphism of $I$. Conversely,
any self--homeomorphism of $I$ is determined by its restriction to $\inte(I)$.
It follows that a choice
of homeomorphism $a:\R \to \inte(I)$ determines an {\em isomorphism\/}:
$$a_*:\homeo(\R) \to \homeo(I)$$

The key observation is that for $\rho:\pi_1(\Sigma) \to \homeo(I)$ a generic representation,
for each $\gamma \in \pi_1(\Sigma)$ and each $p \in \f(\rho(\gamma))$ the germ of
$\rho(\gamma)$ at $p$ is nontrivial.

\begin{lem}\label{generic_exists}
There exist generic representations.
\end{lem}
\begin{proof}
Choose a faithful representation: 
$$\sigma: \pi_1(\Sigma) \to PSL(2,\R)$$
For example, $\sigma$ could be a Fuchsian representation arising from a
hyperbolic structure on $\Sigma$.

$PSL(2,\R)$ can be naturally thought of as a subgroup of $\homeo^+(S^1)$,
the group of orientation--preserving homeomorphisms of the circle,
by considering its action on the circle at infinity of the hyperbolic plane.
For $\gamma \in PSL(2,\R)$ nontrivial, $\gamma$ has $\le 2$ fixed points in
$S^1$. If we choose our representation $\sigma$ corresponding to a hyperbolic structure on
$\Sigma$ where the boundary curves are homotopic to geodesic loops, then
every nontrivial $\gamma$ in the image of $\sigma$ will be {\em hyperbolic}, and
have exactly two fixed points in $S^1$. Moreover, $(\gamma - \id)' \ne 0$ at either
of these two fixed points.

Given a representation $\sigma:\pi_1(\Sigma) \to \homeo^+(S^1)$, we can
construct an oriented {\em foliated} circle bundle $E$ over $\Sigma$ by forming
the quotient of
$$\til{\Sigma} \times S^1$$
foliated by planes $\lbrace \til{\Sigma} \times \text{point} \rbrace$
by the action of $\pi_1(\Sigma)$, where: 
$$\gamma(s,\theta) = (\gamma(s),\sigma(\gamma)(\theta))$$
The Euler class of $\sigma$ is then defined to be the Euler class of this
topological circle bundle, i.e.\ the obstruction to trivializing this
bundle over the $2$-skeleton of $\Sigma$.
Since $\Sigma$ is a punctured surface, it is homotopic to a graph, and
the Euler class of $\sigma$ is $0$. (See \cite{jM58} for a discussion of
Euler classes of representations and circle bundles.)

Thus, $E$ is a {\em trivial} circle bundle $\Sigma \times S^1$, covered by
$\til{E} = \Sigma \times \R$. The foliation of $E$ lifts to a foliation of $\til{E}$
transverse to the $\R$ fibers, and determines a lift of $\sigma$
to a representation:
$$\til{\sigma}: \pi_1(\Sigma) \to \til{PSL(2,\R)}$$

Here $\til{PSL(2,\R)}$ is the universal cover of the group $PSL(2,\R)$. It can be thought
of as a group of homeomorphisms of $\R$ which are periodic with period $2\pi$ and which
cover homeomorphisms of $S^1$ in $PSL(2,\R)$. There is a short exact sequence:
$$0 \to \Z \to \til{PSL(2,\R)} \to PSL(2,\R) \to 0$$

Choose a homeomorphism $a:\R \to I$, and let $a_*\til{\sigma}$ be the
induced representation of $\pi_1(\Sigma)$ in $\homeo(I)$. Then this
representation is clearly generic.
\end{proof}

\begin{lem}\label{typical_construction}
Suppose $\rho:\pi_1(\Sigma) \to \homeo(I)$ is a generic representation.
Then there is a typical foliation $\F_\rho$
of $N$. If $\rho$ is well--tapered, $\F_\rho$ can be chosen to be
$C^\infty$.
\end{lem}
\begin{proof}
The representation $\rho$ gives rise to a foliation of $\Sigma \times I$
by quotienting out the product foliation of $\til{\Sigma} \times I$ by planes
$\lbrace \til{\Sigma} \times \text{point} \rbrace$ by
the action of $\pi_1(\Sigma)$ in the usual way:
$$\gamma: (s,\theta) \to (\gamma(s),\rho(\gamma)(\theta))$$
This is basically the same construction that we used in lemma~\ref{generic_exists}.

Then the top and bottom of $\Sigma \times I$ can be glued up by
the monodromy $\phi$ of $\Sigma$ to give a foliation $\F_\rho$ of $N$.
By construction, $\pi_1$ of the closed leaf $\Sigma \times 0$ injects
into the holonomy groupoid. Moreover, since for any nontrivial
$\alpha$ the element $\rho(\alpha)$ has isolated fixed points in the
interior of $I$, and since every essential closed loop $\alpha'$ in a leaf
$\mu$ of $\F_\rho$ projects to some $\alpha$ by the projection
$\Sigma \times I \to \Sigma$, it follows that the holonomy of $\alpha'$
is nontrivial on either side along $\mu$. Thus $\F_\rho$ is typical.

Clearly, $\F_\rho$ is $C^\infty$ in the interior of $\Sigma \times I$ if
the representation is smooth. Moreover, if $\rho$ is well--tapered, the
holonomy on either side of a loop in the closed surface $\Sigma \times 0$
is infinitely tangent to the identity, and therefore $\F_\rho$ is smooth
there too.
\end{proof}

\begin{thm}
Let $M$ be a closed orientable 
$3$--manifold. Then $M$ admits a typical surface
foliation $\F$.
\end{thm}
\begin{proof}
We find an open--book decomposition of $M$ into $K \coprod N$, by
theorem~\ref{open_book}. We foliate $N$ with a typical foliation $\F_\rho$
by lemma~\ref{generic_exists} and lemma~\ref{typical_construction}.
We twist this foliation around $\partial K$ and fill in $K$ with a
Reeb component to get a foliation $\F$.
Since the meridian of $\partial K$ is transverse to
the circles $\partial \Sigma$, it can be chosen transverse to
$\F_\rho|_{\partial N}$, and therefore the fundamental group of
the torus $\partial K$ injects into the holonomy groupoid of $\F$.
All the leaves in the interior of $K$ are planes, so their fundamental
groups vacuously inject into the holonomy groupoid of $\F$. It follows
that $\F$ is typical.
\end{proof}

\begin{cor}
Any closed orientable $3$--manifold $M$ is a classifying space
$\BG$ for some groupoid $\Gamma$ of germs of homeomorphisms of $\R$.
\end{cor}
\begin{proof}
This follows from Haefliger's theorem~\ref{typical_classifies}.
\end{proof}

\section{Finding well--tapered representations}

In this section we will construct well--tapered representations of
$\pi_1(\Sigma)$ for $\Sigma$ a punctured surface of finite type. In fact, we
will show that well--tapered representations are dense in the space of
representations in $C^\infty$ diffeomorphisms of $I$ infinitely tangent
to the identity at the boundary.

\begin{defn}
Let $C^\infty(I)$ denote the space of $C^\infty$ real--valued functions on $I$.
The {\em $C^r$ topology} on $C^\infty(I)$
is the topology of uniform convergence of the
first $r$ derivatives. The {\em $C^\infty$} topology
is the topology generated by the $C^r$ topology for all $r$. That is,
$f_i \to f$ in the $C^\infty$ topology if the $r$th derivatives $f_i^{(r)} \to
f^{(r)}$ uniformly, for each $r$.
\end{defn}

Note that $C^\infty(I)$ with the $C^\infty$ topology
is {\em not} a Banach space. However, it is a Fr\'echet space --- i.e.\ a
metrizable and complete locally convex topological vector space. (see
\cite{rE65}).

\begin{defn}
Let $\C$ denote the group of $C^\infty$ diffeomorphisms of $I$
infinitely tangent to the identity on $\partial I$.
\end{defn}

The subspace of $C^\infty(I)$ consisting of functions on $I$ infinitely tangent
to the identity on $\partial I$ is an {\em affine} Fr\'echet space.
The group $\C$ can be identified with an open subset of this space
of functions, i.e.\ it is a Fr\'echet manifold. In particular it satisfies
the {\em Baire property} --- any countable intersection of open dense
subsets of $\C$ is dense.

\begin{rmk}
Sergeraert showed in \cite{fS77} that $\C$ is a perfect group.
\end{rmk}

Let $\Sigma$ be a punctured surface of finite type.
Fix a generating set $\alpha_1,\dots,\alpha_n$ for $\pi_1(\Sigma)$. By
hypothesis, the group $\pi_1(\Sigma)$ is free on these generators. The
space of representations $\text{Hom}(\pi_1(\Sigma),\C)$ can therefore
be identitifed with $\C^n$. This space also satisfies the Baire property.

\begin{thm}\label{well_tapered}
Let $\Sigma$ be a punctured surface of finite type.
The set of well--tapered representations $\pi_1(\Sigma) \to \C$ is
dense in the space $\text{Hom}(\pi_1(\Sigma),\C) = \C^n$ of all representations.
\end{thm}
\begin{proof}
The strategy of the proof is as follows. Say an element
$\gamma \in \C$ is {\em robustly generic} if the set of fixed points of
$\gamma$ in $\inte(I)$ has no accumulation points, and if for each fixed point 
$p \in \inte(I)$, $\gamma$ is {\em not} infinitely tangent
to the identity at $p$. For $J \subset \inte(I)$ a {\em closed} subinterval, we say
$\gamma \in \C$ is {\em robustly generic on $J$} if $\f(\gamma) \cap J$ is finite,
and for each $p \in \f(\gamma) \cap J$, $\gamma$ is not infinitely tangent to the
identity at $p$. If 
$$ \cdots \subset J_i \subset J_{i+1} \subset \cdots$$
is an increasing sequence of closed subintervals of $\inte(I)$
whose union is $\inte(I)$, then an element $\gamma \in \C$ is robustly generic iff
it is robustly generic on each $J_i$.

We will show that for each nontrivial $\alpha \in \pi_1(\Sigma)$ and each closed
$J \subset \inte(I)$ the subset of $\rho \in \C^n$ for which $\rho(\alpha)$ is robustly generic
on $J$ is open and dense. If $\rho(\alpha)$ is robustly generic on $J$ for every nontrivial
$\alpha$ and every $J \subset \inte(I)$, then
$\rho$ is a well--tapered representation. It would follow from the Baire property
for $\C^n$ that the set of well--tapered representations is dense --- i.e.\
this would be sufficient to establish the theorem.

Firstly, we establish openness of the condition. Let $J \subset \inte(I)$ be closed, and
suppose $\gamma \in \C$ is robustly generic on $J$. Then $\f(\gamma) \cap J$ is a finite
set. For any $\epsilon>0$, if $\nu$ is sufficiently close to $\gamma$,
we can estimate that
$\f(\nu) \cap J \subset N_\epsilon(\f(\gamma)) \cap J$, where $N_\epsilon(\cdot)$
denotes the $\epsilon$--neighborhood of a set. Moreover, we can require
that $\nu^{(n)}$ is arbitrarily close to $\gamma^{(n)}$ on this neighborhood, for any
fixed $n$. Notice that the fixed points of $\nu$ are exactly the
zeroes of $\nu - \text{id}$. By the robustness hypothesis,
for each $p \in \f(\nu) \cap J$, and for some $\delta>0$, there is
a $\kappa>0$ such that for some smallest $n$,
$$|(\nu-\text{id})^{(n)}(q)|>\delta$$ for all $q \in N_\kappa(p)$, and
$$(\nu-\text{id})^{(m)}(p)=0$$ for all $m<n$. We change co--ordinates near $p$
by setting $y= (x-p)^n$, and see that with respect to $y$:
$$\Biggl| \frac d {dy} (\nu - \text{id})(0) \Biggr| > 0$$
In particular, $p$ is an {\em isolated} point of $\f(\nu)$, and therefore
$\nu$ is robustly generic.

It remains to establish density. Pick a closed interval $J \subset\inte(I)$ and
some $\alpha \in \pi_1(\Sigma)$, and
suppose we have chosen some $\rho \in \C^n$ for which $\rho(\alpha)$ is not robustly 
generic on $J$. Fix a degree of smoothness $m$. For each $\epsilon>0$,
we can find a collection of closed intervals $J_i \subset \inte(I)$
for $1\le i\le n$ such that for
$\sigma \in \C^n$ chosen $\epsilon$--close to $\rho$ in the $C^m$ norm,
the restriction of $\sigma(\alpha)$ to $J$ is determined by the
restrictions of $\sigma(\alpha_i)$ to $J_i$. Moreover, for each
$\delta$, we can find an $\epsilon$ such that if $\sigma$ is $\epsilon$--close
to $\rho$ on the generating set $\lbrace \alpha_i \rbrace$, it is $\delta$ close to $\rho$ on
$\alpha$, again in the $C^m$ norm.

We will show how to modify $\rho(\alpha_i)$
on the subintervals $J_i$ so that the modified $\rho(\alpha)$ is
robustly generic on $J$. Think of the image of each generator
$\alpha_i$ under $\rho$ as a function from $I$ to $I$:
$$\rho(\alpha_i)(x) = f_i(x)$$
Fix a positive integer $N\gg m$. Let $X_N(i)$ denote the space of
invertible $C^m$ functions $g_i:I \to I$ with the following properties:
\begin{enumerate}
\item{The function $g_i$ agrees with $f_i$ outside $J_i$.}
\item{$g_i$ is tangent to order $m$ to $f_i$ at $\partial J_i$.}
\item{The restriction $g_i|_{J_i}$ is a real polynomial of degree $\le N$.} 
\end{enumerate}
The $X_N(i)$ determine a subspace of $\C^n$, which we denote by $X_N$. Clearly,
$X_N$ is a finite dimensional real analytic variety.

For any $\sigma \in X_N$,
the restriction $\sigma(\alpha)|_J$ is a real analytic function, which
extends to a holomorphic (complex--valued) function in a neighborhood of $J \subset
\mathbb{C}$. In particular, if $\sigma(\alpha)$ is nontrivial on $J$,
it is robustly generic there. Let $Y_N \subset X_N$ be the
analytic subvariety defined by the condition that $\sigma(\alpha)|_J$ is trivial iff
$\sigma \in Y_N$. For sufficiently large $N$, we claim $X_N\backslash Y_N$
is nonempty, and therefore is
open and dense in $X_N$. Let $\tau:\pi_1(\Sigma) \to \homeo(I)$ denote
a generic representation constructed in lemma~\ref{generic_exists}. $\tau$ can
be chosen to be real analytic on $\inte(I)$, by choosing some real analytic homeomorphism
$a:\R \to \inte(I)$. We choose closed subintervals $K_i \subset \inte(J_i)$
and let $\sigma(\alpha_i)$ be $C^m$ close to $\tau(\alpha_i)$ on $K_i$. If the
$K_i$ are large enough, they determine the value of $\sigma(\alpha)$ on some nonempty
$K \subset J$. Since $\tau(\alpha)|_K$ is nontrivial (i.e.\ not equal to the identity there),
it follows that $\sigma(\alpha)|_K$ is nontrivial. So $X_N\backslash Y_N$ is
nonempty for large $N$ and a suitable choice of $J_i$, and therefore it is open
and dense in $X_N$. It follows that we can find $\sigma \in X_N$ which is $C^m$ close
to $\rho$ for which $\sigma(\alpha)|_J$ is robustly generic.

If we $C^N$ approximate $\sigma$ by $\sigma'\in \C^n$ 
to make it $C^\infty$, $\sigma'(\alpha)|_J$ is still robustly generic,
by openness. It follows that we have established that the set of $\sigma \in \C^n$ for
which $\sigma(\alpha)$ is robustly generic on $J$ is dense.

Let $J_i$ be an increasing sequence of closed subsets of $\inte(I)$ whose union is
$\inte(I)$. Enumerate the nontrivial $\alpha \in \pi_1(\Sigma)$ in some order $\alpha_i$,
where we can assume if we like that this ordering agrees with the labels on the generating set.
For each $i$, the set of $\sigma \in \C^n$ for which $\alpha_j$ is robustly generic
on $J_i$ for all $j \le i$ is open and dense. The intersection of this countable family
of open, dense sets is dense,
and therefore there are a dense set of $\sigma$ for which
$\sigma(\alpha)$ is robustly generic for {\em every} nontrivial $\alpha$.
Such a $\sigma$ is well--tapered.
\end{proof}

Notice that the proof actually tells us slightly more than the conclusion of the
theorem. Given $\sigma \in \C^n$, we can find
$\sigma' \in \C^n$ robustly generic which is obtained by successively perturbing
$\sigma|_{J_i}$ on intervals $J_i \to \inte(I)$
by smaller and smaller amounts. So for any $\alpha \in \pi_1(\Sigma)$
we can choose these perturbations so that $\sigma'(\alpha) - \sigma(\alpha)$
tapers off as fast as desired as it approaches $\partial I$. Here our notation
suggests that we are thinking of $\sigma(\alpha)$ and $\sigma'(\alpha)$ as functions
from $I$ to $I$. In particular, if $\sigma(\alpha)$ satisfies certain desirable properties ---
for instance, if $\inte(I) \cap \f(\sigma(\alpha))$ is empty --- then this property may be
inherited by $\sigma'$ if we desire.

\begin{defn}
A function $f:I \to \R$ is {\em stretchable} if
it satisfies the following conditions:
\begin{enumerate}
\item{$f$ is $C^\infty$ and $f$ is positive on $\inte(I)$.}
\item{$f$ is infinitely tangent to $0$ at the ends of $I$.}

\noindent{Furthermore, let $X$ denote the vector field $f(t)\frac d {dt}$ on $I$, and
$\phi_t:I \to I$ denote the diffeomorphism generated by time $t$
flow along $X$. Define $\alpha:I \to I$ by the properties that
$\alpha(\frac 1 2) = \frac 1 2$ and
$\alpha \phi_{2t} \alpha^{-1} = \phi_t$. Then}
\item{$\alpha \in \C$.}
\end{enumerate}
We say in this case that $\phi_t$ is stretchable.
\end{defn}

Stretchable functions certainly exist. By basic existence and uniqueness
results for differential equations, $\alpha$ exists and is $C^\infty$ on
the interior of $I$; in order to ensure the correct boundary behavior
of $\alpha$, we just need $f$ to taper off sufficiently fast as it approaches
$\partial I$.

\begin{cor}\label{damped_ends}
Let $\Sigma$ be a hyperbolic punctured surface of finite type.
There is a well--tapered representation $\rho:\pi_1(\Sigma) \to \C$
such that if $\gamma$ denotes
the image in $\pi_1(\Sigma)$ of a small loop around a puncture, the
element $\rho(\gamma)$ fixes no points in the interior of $I$, and
$\rho(\gamma) - \text{id}$ tapers off as fast as desired as it approaches
$\partial I$.
\end{cor}
\begin{proof}
Sergeraert's theorem implies that any element of $\C$ can be written
as a product of commutators; however, this observation does not help us,
since we have no bound on the number of commutators needed to express
a given element of $\C$.

On the other hand, for $f,X,\phi_{2t},\alpha$ as above with $\phi_t$
stretchable, the commutator $$[\phi_{2t},\alpha] = \phi_t$$ fixes no points
in the interior of $I$. Since $\gamma$ is a product of commutators,
we can choose a representation $\sigma$ taking the first commutator to
$\phi_t$ and all the others to $\text{id}$.
Then any $\rho$ sufficiently close to $\sigma$ will satisfy
the conclustion of the corollary; such well--tapered $\rho$ exist by
theorem~\ref{well_tapered}.
\end{proof}

\section{Pixton actions}

\begin{thm}
Let $M$ be an orientable $3$--manifold. Then $M$ admits a $C^1$
typical surface foliation $\F$.
\end{thm}
\begin{proof}
As before, we decompose $M = K \coprod N$. We foliate $N$ by using a
well--tapered representation of the fundamental group of a fiber,
by theorem~\ref{well_tapered}. Since the dynamics of the
holonomy of $\partial \Sigma$ has no fixed points in the interior of $I$,
the end of $N$ has a single cylinder leaf, and plane leaves spiralling
towards it, in a $C^\infty$ fashion. This can be twisted around $\partial K$
while keeping it $C^1$, where it will be eventually tangent to
the suspension of a Pixton representation of $\Z\oplus \Z$ in the
group of $C^1$ homeomorphisms of $I$; see \cite{Pix}. With notation
from corollary~\ref{damped_ends}, as we spin
the end of $N$ around $\partial K$, the holonomy is generated by
some $\beta:I \to I$ with no fixed points near $0$ and
the suspension $\alpha$ of $\rho(\gamma)$ by $\beta$. Explicitly: 
there is an inclusion $\varphi:I \to [p_1,p_0]$ where $0<p_1 < p_0 < 1$. 
Choose a sequence $p_i \to 0$ where $\beta(p_i) = p_{i+1}$ and define $\alpha'$
by $$\alpha'|_{[p_1,p_0]} = \varphi \rho(\gamma) \varphi^{-1}$$
and $\alpha' = \text{id}$ outside $[p_1,p_0]$.
Then define: $$\alpha = \prod_{i=-\infty}^\infty \beta^i \alpha' \beta^{-i}$$
Such an $\alpha$ is the suspension of $\rho(\gamma)$ by $\beta$; clearly,
the elements $\alpha$ and $\beta$ commute. We would like to pick $\beta$
and $p_i$ such that $\alpha$ and $\beta$ are both tangent to first order to
the identity at $0$. First of all, we should choose the $p_i$ so that
the ratios $[p_i,p_{i+1}]/[p_{i-1},p_i] \to 1$ as $i \to \infty$, for instance
by choosing $p_i = \frac 1 i$. The restriction of $\alpha$ to 
$(p_i,p_{i+1})$ is conjugate to a translation of $\R$; so we just need to
choose a sequence of translations on each of the $(p_i,p_{i+1})$ which,
after rescaling, converge geometrically in the $C^1$ topology to the identity.
This requires $\beta$ to expand the intervals where $(\alpha-\text{id})'$ 
is very small, and contract the intervals where 
$(\alpha-\text{id})'$ is larger. Let $\|f\|_\infty^i$ denote the $L^\infty$
norm of $f$ on the interval $[p_i,p_{i+1}]$.
If $\rho(\gamma)$ decays sufficiently quickly
to the identity near $\partial I$, we can find, for any $\delta_i>0$,
a restriction of $\beta$ to $[p_i,p_{i+1}]$ for which
$\|\beta'\|_\infty^i <\delta_i$ but which satisfies
$$\frac {\|(\alpha -\text{id})'\|_\infty^{i+1}} 
{\|(\alpha -\text{id})'\|_\infty^i} \le 1 -\epsilon_i$$
for some definite $\epsilon_i > 0$.
Choose a sequence $\delta_i \to 0$ so that
$\prod_i (1-\epsilon_i) \to 0$.
More details are found in \cite{Pix}.
\end{proof}

On the other hand, this foliation cannot be made $C^2$, by Kopell's lemma
(\cite{Kop}), since the holonomy representation of
$\partial \Sigma$ on $\partial N$, thought of as a circle bundle over
$\partial \Sigma$, has a fixed point. Thus, although the foliation of
$N$ is $C^\infty$, this foliation cannot be smoothly twisted around the
boundary torus and plugged in with a Reeb component.

Tantalisingly, the foliation $\F$ of $N$ {\em can} be modified slightly
so that it can be smoothly spun around $\partial N$:

\begin{thm}
The foliation $\F$ of $N$ as above can be $C^\infty$ approximated by
a foliation $\G$ such that the holonomy of $\partial \Sigma$ has no
fixed points.
\end{thm}
\begin{proof}
Let $\gamma$ be an essential simple arc in $\Sigma$. We cut open the
foliation $\F$ along a small transverse foliated rectangle $R$,
shift slightly and reglue. Since the holonomy along $\gamma$ followed
by an arc in $\partial \Sigma$ is nontrivial, this perturbation can
be made to modify the holonomy of $\partial \Sigma$ to have no fixed
points. Herman showed in \cite{Her} that for a subset $A \subset S^1$ of
full measure, a $C^\infty$ diffeomorphism from $S^1$ to itself with
rotation number in $A$ is {\em $C^\infty$ conjugate} to a rotation.
It follows that for suitable choice of generic shift along $R$, the rotation
number of this boundary holonomy will lie in $A$, and therefore
this element will be $C^\infty$ conjugate to a rotation.
Such a foliation can be spun {\em smoothly} around a Reeb component
to give a $C^\infty$ foliation of $M$.
\end{proof}

One might ask whether in this context, a ``generic'' perturbation of
the resulting foliation, preserving the fixed--point free irrational
dynamics on the boundary, can make the foliation typical, following
the proof of theorem~\ref{well_tapered}. Unfortunately,
after shearing along $R$, the foliation of $N$ is no longer necessarily
complete; in fact, it will almost certainly blow up holonomy of some
intervals in finite time. This
is not to say that no such perturbation could exist, but new techniques
would be required to find it.

\Addressesr

\end{document}